\documentclass[a4paper,12pt]{article}
\usepackage{amssymb}
\usepackage{amsmath}
\usepackage{amscd}
\usepackage[arrow,matrix]{xy}
\usepackage[bookmarksnumbered,colorlinks,linkcolor=black,citecolor=black,urlcolor=black]{hyperref}
\usepackage{mathtools}

\newtheorem{theo}{Theorem}[section]
\newtheorem{prop}[theo]{Proposition}
\newtheorem{lem}[theo]{Lemma}
\newtheorem{cor}[theo]{Corollary}
\newtheorem{defi}[theo]{Definition}
\newtheorem{rem}[theo]{Remark}
\newtheorem{exe}[theo]{Exercise}
\newtheorem{exa}[theo]{Example}

\newtheorem{conj}[theo]{Conjecture}
\newtheorem{ques}[theo]{Question}

\newcommand{\bthe}{\begin{theo}}
\newcommand{\ble}{\begin{lem}}
\newcommand{\bpr}{\begin{prop}}
\newcommand{\bco}{\begin{cor}}
\newcommand{\bde}{\begin{defi}}
\newcommand{\ethe}{\end{theo}}
\newcommand{\ele}{\end{lem}}
\newcommand{\epr}{\end{prop}}
\newcommand{\eco}{\end{cor}}
\newcommand{\ede}{\end{defi}}
\newcommand{\brem}{\begin{rem}}
\newcommand{\erem}{\end{rem}}

\newcommand{\bexe}{\begin{exe}}
\newcommand{\eexe}{\end{exe}}

\newcommand{\bexa}{\begin{exa}}
\newcommand{\eexa}{\end{exa}}

\newcommand{\bconj}{\begin{conj}}
\newcommand{\econj}{\end{conj}}
\newcommand{\bques}{\begin{ques}}
\newcommand{\eques}{\end{ques}}

\def \Br {{\rm{Br}}}

\def \si {{\sigma}}

\def \A{{\mathbb A}}
\def \P{{\mathbb P}}

\def \Spec {{\rm{Spec}}}
\def \dim {{\rm{dim}}}

\def \GL {{\rm {GL}}}
\def \SL {{\rm {SL}}}

\def \PGL {{\rm {PGL}}}

\def \tr {{\rm {tr}}}

\def \Z {{\mathbb Z}}

\def \H {{\rm H}}
\def\G{{\mathbb G}}

\def\T{{\cal T}}

\def\tr{{\rm tr}}

\def\sA{{\cal A}}

\def\sF{{\cal F}}

\def\si{\sigma}

\def\et{{\rm{\acute et}}}

\DeclareFontEncoding{OT2}{}{} 
\DeclareTextFontCommand{\textcyr}{\fontencoding{OT2}\fontfamily{wncyr}\fontseries{m}\fontshape{n}\selectfont}

\usepackage[textheight=220mm,textwidth=150mm]{geometry}

\title{Birational properties of word varieties}
\author{Tatiana Bandman, Boris Kunyavski\u\i, and Alexei N. Skorobogatov}

\date{\today}
\begin{document}
\maketitle

\centerline{\em to Yuri Tschinkel, on his 60-th birthday}

\begin{abstract}
\noindent We prove that the subvariety of $\SL(2)\times\SL(2)$ given by the
matrix equation $w(X,Y)=\alpha$, where $w$ is a word in two letters,
is closely related to an explicit smooth conic bundle over the associated 
`trace surface' in the 3-dimensional affine space. When $w$ is the commutator word,
we show that this variety can be irrational if the ground field $k$ is not algebraically closed,
answering a question of Rapinchuk, Benyash-Krivetz, and Chernousov.
When $k$ is a number field,
it satisfies weak approximation with the Brauer--Manin obstruction.
\end{abstract}

\section*{Introduction}

Let $w(x,y)$ be a non-trivial word in two letters $x$ and $y$, that is, an element of the free
group ${\mathcal F}_2$ with two generators.
Let $k$ be a field of characteristic different from 2.
For $\alpha\in\SL(2)(k)$ we define the {\em word variety}
$S_{w,\alpha}$ as the closed subvariety of $\SL(2)\times\SL(2)$ given by the matrix equation $w(A,B)=\alpha$. For a survey of the extensive literature on word equations, see \cite{GKP}.
The motivating question of this paper concerns birational properties of $S_{w,\alpha}$.
We also obtain results about $k$-rational points on $S_{w,\alpha}$ when $k$ is a number field.

One much studied example is the commutator
variety $ABA^{-1}B^{-1}=\alpha$, see \cite{Th}, \cite{RBKC} and, 
more recently, \cite{LL21}, \cite{GS22, GMS21}. 
When $\alpha$ is semisimple and non-central,
we show that the commutator variety is a dense open subset of a smooth conic bundle
(a Severi--Brauer scheme of relative dimension 1)
over the Markoff surface in $\A^3_k$ with equation
\begin{equation}
s^2+t^2+u^2-stu=\tr(\alpha)+2.\label{M}
\end{equation}
The associated Brauer class is $(s^2-4,\tr(\alpha)-2)$; it
generates the unramified Brauer group of the
`generic' Markoff surface modulo $\Br(k)$. 
In particular, the commutator variety
is birationally equivalent to the affine subvariety of $\A^5_k$ given by
$$s^2+t^2+u^2-stu-(\tr(\alpha)+2)=s^2-4 +(\tr(\alpha)-2)x^2-y^2=0.$$
We prove that the commutator variety is $k$-rational if and only if
the Markoff surface (\ref{M}) is $k$-rational, which is the case if and only if either $\tr(\alpha)-2$ or 
$\tr(\alpha)^2-4$ is a square in $k$, see Theorem \ref{3.3}.
In particular, for `general' $\alpha\in \SL(2)(k)$  the commutator variety is irrational over $k$. 
This gives a negative answer to a question of 
Rapinchuk, Benyash-Krivetz, and Chernousov \cite[page 50]{RBKC}. 
When $k$ is a number field, we show that 
the Brauer--Manin obstruction is the only obstruction to weak approximation 
on smooth and proper models of the commutator variety, see Proposition \ref{BMO}.

Our method is based on the systematic use of the 
`trace polynomial' of the word $w(x,y)$. This is very classical and goes back to Vogt, Klein, and Fricke.
Using the Cayley--Hamilton theorem it is easy to see that the function
$\SL(2)\times\SL(2)\to \A^1_k$ given by $\tr(w(A,B))$ is a polynomial $P_w(s,t,u)$ in
the variables
$s=\tr(A)$, $t=\tr(B)$, $u=\tr(AB)$. We define the {\em trace surface} $H_{w,a}$, for $a\in k$,
as the surface in $\A^3_k$ given by $P_w(s,t,u)=a$. If $w(x,y)$ belongs to the commutator
subgroup of the free group on two letters $x$ and $y$, and $\alpha\in\SL(2)(k)$ is semisimple and non-central
with $\tr(\alpha)=a$, the word variety $S_{w,\alpha}$ is a torsor over the trace surface $H_{w,a}$
for the norm 1 torus associated to 
the splitting field $k(\sqrt{a^2-4})$ of the characteristic polynomial of $\alpha$.
This implies that $S_{w,\alpha}$ is a dense open subset of an
explicit smooth conic bundle over $H_{w,a}$, see Theorem \ref{t1}. 
This largely reduces the study of solutions of the word equation $w(A,B)=\alpha$ to the study
of rational points on the affine surface $H_{w,a}$.

We are grateful to Tim Browning, Jean-Louis Colliot-Th\'el\`ene, David
Harari, and Yuri Zarhin for helpful discussions and encouragement.
We thank the anonymous referee for their careful reading of the paper.
The research of the second named author was supported by the Israel Science
Foundation grant 1994/20. The work on this paper benefited from the authors' multiple
visits to the Max Planck Institute for Mathematics in Bonn and 
the last named author's stay at Pennsylvania State University as a Shapiro Visitor.
The hospitality and support of these institutions are gratefully acknowledged.

\section{Preliminaries}

\subsection{Reduction of the structure group} \label{red}

For the reader's convenience we recall the definition of the push-forward
of torsors, see \cite[\S 2.2]{S01} and \cite[Proposition III.3.2.1]{Giraud}.
Let $G_1$ be a closed subgroup of an affine algebraic group 
$G$ over a field $k$.
Let $Z$ be a variety over $k$ and let $\T_1\to Z$ be a right $Z$-torsor for $G_1$ for the \'etale topology.
The class of this torsor is an element $[\T_1]$ of the \'etale cohomology set $\H^1_\et(Z,G_1)$.
The {\em push-forward} of $\T_1$ along the morphism $G_1\hookrightarrow G$
is the quotient $\T$ of $\T_1\times_k G$ by the diagonal action of $G_1$,
where $G_1$ acts on $\T_1$ on the right and on $G$ on the left.
Then $\T\to Z$ inherits a right action of $G$ making it a $Z$-torsor for $G$.
This gives a map of pointed sets $\H^1_\et(Z,G_1)\to \H^1_\et(Z,G)$
that sends $[\T_1]$ to $[\T]$. When $Z$ is $\Spec(k)$, the \'etale cohomology sets
become Galois cohomology sets. 

If a $Z$-torsor $\T\to Z$ for $G$ is the push-forward
of some $Z$-torsor for $G_1$, then one says that $\T$ {\em lifts} to a 
torsor for $G_1$, or that the structure
group of $\T$ {\em reduces} to $G_1$. 

If $G_1$ is commutative, every right $Z$-torsor 
for $G_1$ is also a left $Z$-torsor for $G_1$. (In general, a right torsor for $G_1$ is a left torsor for 
an inner form of $G_1$.) In this case we have the following criterion for the reduction of the structure
group from $G$ to $G_1$.

\bpr \label{comm}
Let $G_1$ be a closed subgroup of an affine algebraic group $G$ defined over a field $k$. 
Assume that $G_1$ is commutative. Let $Z$ be a variety over $k$
and let $\T$ be a right $Z$-torsor for $G$. The structure group of $\T$ reduces 
to $G_1$ if and only if there is a $G$-equivariant morphism of $Z$-schemes 
$\T\to (G_1\backslash G)\times_k Z$.
\epr
{\em Proof.} Let $\T_1$ be a right $Z$-torsor for $G_1$. Let $\T=(\T_1\times_k G)/G_1$
be the quotient by the diagonal action of $G_1$,
which acts on $\T_1$ on the right and on $G$ on the left. Thus $\T$ has a right action
of $G$ making it a $Z$-torsor for $G$, but also a left action of $G_1$. Let $G_1\backslash\T$
be the quotient by this left action of $G_1$. We have isomorphisms
$$G_1\backslash\T\cong G_1\backslash (\T_1\times_k G)/G_1\cong Z\times_k (G_1\backslash G),$$
since the two actions of $G_1$ obviously commute, and $G_1\backslash \T_1=Z$.

Conversely, let $\varphi\colon \T\to (G_1\backslash G)\times_k Z$ be a 
$G$-equivariant morphism of $Z$-schemes, with $G$ acting on the right. 
Let $x_0$ be the $k$-point of 
$G_1\backslash G$ given by the identity element of $G(k)$. Note that the right action of $G_1$ on 
$G_1\backslash G$ preserves $x_0$.
Thus $\T_1:=\varphi^{-1}(x_0\times_k Z)$ is a closed subvariety of $\T$ stable under the
right action of $G_1$ so that $\T_1\hookrightarrow \T$ is compatible with the injective
homomorphism $G_1\hookrightarrow G$. Since $\T$ is a right $Z$-torsor for $G$,
this implies that $\T_1$ is a right $Z$-torsor for $G_1$
such that $\T$ is the push-forward of $\T_1$ along $G_1\hookrightarrow G$. \hfill $\Box$

\subsection{Norm $1$ tori for quadratic extensions} \label{1.1}

Let $k$ be a field of characteristic different from 2. For a quadratic extension $L/k$ we write
$R_{L/k}(\G_{m,L})$ for the Weil restriction of $\G_{m,L}$.
The attached norm 1 torus $R^1_{L/k}(\G_{m,L})$ is defined by the exact sequence of $k$-tori
\begin{equation}
1\to R^1_{L/k}(\G_{m,L})\to R_{L/k}(\G_{m,L})\xrightarrow{\rm N} \G_{m,k}\to 1,
\label{norm}\end{equation}
where N is induced by the norm map $L\to k$. As a variety, $R^1_{L/k}(\G_{m,L})$
is the affine conic $x^2-ay^2=1$, where $a\in k^\times$ is such that $L=k(\sqrt{a})$.
The long exact sequence 
of Galois cohomology groups attached to (\ref{norm}), by Hilbert's Theorem 90, gives rise to
a canonical isomorphism
$$\H^1(k,R^1_{L/k}(\G_{m,L}))\cong k^\times/{\rm N}(L^\times).$$
This group classifies isomorphism classes of $k$-torsors for $R^1_{L/k}(\G_{m,L})$.
Such a torsor is an affine conic $x^2-ay^2=b$, where $b\in k^\times$, with the
obvious action of $R^1_{L/k}(\G_{m,L})$. Its class
in $\H^1(k,R^1_{L/k}(\G_{m,L}))$ is the image of $b$ in $k^\times/{\rm N}(L^\times)$.


\subsection{Affine and projective conics} \label{1.3}

Let $Z$ be a $k$-variety and let $\T\to Z$ be a right $Z$-torsor for $\PGL(n)$.
To the class $[\T]\in\H^1_\et(Z,\PGL(n))$ one canonically associates an element $\partial([\T])$
of the (cohomological) Brauer group
$\Br(Z)=\H^2_\et(Z,\G_m)$, where $\partial$ is the connecting map
attached to the central extension of algebraic groups
$$1\to\G_m\to\GL(n)\to\PGL(n)\to 1. $$
To the torsor $\T\to Z$ one associates a Severi--Brauer $Z$-scheme of relative dimension $n-1$,
as follows. The standard $n$-dimensional representation of $\GL(n)$
gives rise to a transitive action of $\PGL(n)$ on the projective space $\P^{n-1}_k$, which we
write as a left action.
Let $X$ be the quotient of $\T\times_k\P^{n-1}_k$ by the diagonal action of $\PGL(n)$,
acting on $\T$ on the right and on $\P^{n-1}_k$ on the left.
The $Z$-scheme $X$ is \'etale locally isomorphic to $Z\times_k\P^{n-1}_k$, so $X\to Z$ is a 
Severi--Brauer scheme of relative dimension $n-1$.
The class of $X$ in $\Br(Z)$ is defined as $[X]=\partial([\T])\in\Br(Z)$.

\ble \label{lem4}
Let $Z$ be a variety over $k$, let $\T\to Z$ be a right $Z$-torsor for $\PGL(n)$, and let
$X=(\T\times_k\P^{n-1}_k)/\PGL(n)$ be the attached Severi--Brauer $Z$-scheme  
of relative dimension $n-1$. If
$\T$ lifts to a $Z$-torsor $\T_1$ for a maximal $k$-torus $T\subset \PGL(n)$, then
there is an open embedding of $Z$-schemes $\T_1\hookrightarrow X$ with dense image.
\ele
{\em Proof.} Since $\T$ is the push-forward of $\T_1$, we have canonical isomorphisms
$$X=(\T\times_k\P^{n-1}_k)/\PGL(n)\cong 
(\T_1\times_k \PGL(n)\times_k\P^{n-1}_k)/(T\times_k\PGL(n))\cong
(\T_1\times_k\P^{n-1}_k)/T,$$
where $T$ acts on $\P^{n-1}_k$ on the left as a subgroup of $\PGL(n)$.
Indeed, the actions of $T$ and $\PGL(n)$ on $\T_1\times_k \PGL(n)\times_k\P^{n-1}_k$ commute,
as immediately follows from their definitions.

The restriction of the action of $\PGL(n)$ on $\P^{n-1}_k$ to $T$
has a dense orbit on which $T$ acts freely.
Fixing a $k$-point in this orbit gives an open embedding $T\hookrightarrow \P^{n-1}_k$.
We identify $T$ with the image of this embedding, so that $T$ acts on itself by translations.
Twisted by $\T_1$, the embedding $T\hookrightarrow \P^{n-1}_k$ gives rise to the desired
open embedding of $Z$-schemes $\T_1\hookrightarrow X$. \hfill$\Box$

\bexa{\rm \label{dim2}
For $n=2$ and $Z=\Spec(k)$
we can make Lemma \ref{lem4} explicit. Let $T\subset\PGL(2)$ be a maximal $k$-torus
and let $T'$ be its preimage in $\SL(2)$. Since the centre 
$\mu_2$ of $\SL(2)$ is the 2-torsion subgroup $T'[2]\subset T'$, the multiplication by 2 map
gives an isomorphism $T'\cong T'/\mu_2\cong T$. 
So the maximal $k$-tori in $\PGL(2)$ are also the maximal $k$-tori in $\SL(2)$.

For an element $g\in\SL(2)(k)$ we write
$C_{\SL(2)}(g)$ for the centraliser of $g$ in $\SL(2)$.
The torus $T\subset\SL(2)$ is the centraliser $C_{\SL(2)}(g)$ of any non-central $g\in T(k)$.
For a non-central semisimple element $g\in\SL(2)(k)$ with trace $t=\tr(g)$ we have
$$T=C_{\SL(2)}(g)\cong R^1_{L/k}(\G_{m,L}),$$
when $L=k(\sqrt{t^2-4})$ is a quadratic extension of $k$. 
When $t^2-4$ is a square in $k^\times$, we define $L=k\oplus k$. In this case 
$T=R^1_{L/k}(\G_{m,L})$ is the split torus $\G_{m,k}$. It is easy to see that every
field extension of $k$ of degree at most 2 arises in this way.

Let $L=k(\sqrt{a})$ where $a\in k^\times$.
A $k$-torsor $\T_1$ for $R^1_{L/k}(\G_{m,L})$ can be given by
$$x^2-ay^2=b, $$
for some $b\in k^\times$. Let $\T$ be a $k$-torsor for $\PGL(2)$ which is the push-forward of $\T_1$
along any embedding $R^1_{L/k}(\G_{m,L})\hookrightarrow\PGL(2)$, and let
$X$ be the Severi--Brauer
variety associated to $\T$. By Lemma \ref{lem4}, $\T_1$ is a dense open subset of $X$, so $X$
is the projective conic given by the homogenised equation 
$$x^2-ay^2=bz^2. $$
We note that $\T_1$ is the complement to the closed point $\Spec(L)\subset X$ given by $z=0$.
}
\eexa

\brem{\rm \label{uni}
In the above construction we can replace a maximal torus in $\PGL(2)$ by the centraliser
of an element $g\in \SL(2)(k)$ 
such that $\tr(g)=\pm 2$ and $g\neq \pm I$. In this case $G:=C_{\PGL(2)}(g)$
is isomorphic to the additive group $\G_{a,k}$. 
The restriction of the action of $\PGL(2)$ on $\P^1_k$ to $G$
has a dense orbit isomorphic to $\A^1_k$ on which $G$ acts freely. Thus we have an analogue of
Lemma \ref{lem4} in the unipotent case (with the same notation):
if $\T$ lifts to a $Z$-torsor $\T_1$ for $G$, then
there is an open embedding $\T_1\hookrightarrow X$. 
When $Z=\Spec(k)$ we have $\H^1(k,\G_{a,k})=0$ by the additive version of Hilbert's Theorem 90,
hence any $k$-torsor for $G$ is trivial
and so is isomorphic to $\A^1_k$. This implies that $X\cong\P^1_k$.
}\erem

\section{Simultaneous similarity of two matrices}

\subsection{Markoff surfaces and torsors for $\PGL(2)$}

Let $k$ be a field of characteristic different from 2 with algebraic closure $\bar k$. 
Write
$$F(s,t,u)=s^2+t^2+u^2-stu-4.$$
For $d\in k$ let $M_d\subset\A^3_k$ be the affine cubic surface given by $F(s,t,u)=d$,
called a {\em Markoff surface}.
The surface $M_d$ is smooth if and only if $d\neq 0$ and $d\neq -4$.
The singular cubic surface $M_0$, also called the {\em Cayley cubic}, 
has four singular points with coordinates
$s,t,u=\pm 2$ with the product of signs equal to 1.

The Cayley cubic $M_0$ naturally arises in the problem of 
classification of pairs of $(2\times 2)$-matrices up to simultaneous similarity. Write
$$f\colon \SL(2)\times_k\SL(2)\to \A^3_k$$
for the morphism sending $A, B\in \SL(2)(\bar k)$ to $(s,t,u)\in\A^3_k(\bar k)$, where
$$s=\tr(A), \quad t=\tr(B), \quad u=\tr(AB).$$
The action of $\PGL(2)$ by simultaneous conjugation preserves the fibres of $f$.

The following proposition is well-known.
\bpr \label{prop}
Let $A, B\in \SL(2)(\bar k)$. The following properties are equivalent:

{\rm (1)} $A$ and $B$ have a common eigenvector;

{\rm (2)} $\det(AB-BA)=0$;

{\rm (3)} $\tr(ABA^{-1}B^{-1})=2$;

{\rm (4)} $(s,t,u)$ is a $\bar k$-point of the Cayley cubic $M_0$.

\noindent These conditions are satisfied when the centralisers of $A$ and $B$ in $\PGL(2)$
have non-trivial intersection.
\epr
{\em Proof.} The equivalence of (1) and (2) is \cite[Theorem 3.1]{Shemesh}. The equivalence 
of (2) and (3) is elementary, as both statements are equivalent to 
1 being a root of the characteristic polynomial of $ABA^{-1}B^{-1}$.
The equivalence of (3) and (4) follows from the Cayley--Hamilton
theorem 
\begin{equation}
A^2=sA-I, \quad B^2=tB-I, \quad ABAB=uAB-I, \label{CH}
\end{equation}
from which one obtains the Fricke identity 
\begin{equation}
\tr(ABA^{-1}B^{-1})=s^2+t^2+u^2-stu-2. \label{Fricke}
\end{equation}
See also \cite[Theorem 2.9]{Fri83}. For the last statement, let $g\neq \pm I$ be an element
of $\SL(2)(\bar k)$ commuting with $A$ and $B$. The centraliser of $g$ in $\SL(2)(\bar k)$
is commutative, hence $AB=BA$, so (2) holds. \hfill$\Box$

\medskip

Let $V=\A^3_k\setminus M_0$.
If we identify the $k$-points of $\SL(2)\times_k\SL(2)$ with 
representations of the free group ${\mathcal F}_2\to\SL(2)(k)$, then
Proposition \ref{prop} implies that the open subscheme 
$$\T:=f^{-1}(V)\subset \SL(2)\times_k\SL(2)$$ 
parameterises absolutely irreducible representations of ${\mathcal F}_2$, that is, those
that remain irreducible over $\bar k$.

\ble \label{2.1}
The morphism $\T\to V$ is a torsor for $\PGL(2)$ for the \'etale topology.
\ele
{\em Proof.} This is proved in  \cite[Corollaries 6.5, 6.8]{N}. \hfill $\Box$

\medskip

The class of $\T$ is an element $[\T]$ of the \'etale cohomology set $\H^1_\et(V,\PGL(2))$.

\subsection{Affine and projective conic bundles} \label{2.2}

As in \S \ref{1.3}, we associate to the torsor $\T\to V$ a 
Severi--Brauer scheme $X\to V$ of relative dimension 1, that is, a smooth
conic bundle. We would like to compute $[X]\in\Br(V)$.

Let $K=k(s,t,u)$ be the field of functions on $\A^3_k$. Because of the canonical embedding
$\Br(V)\subset\Br(K)$, see \cite[Theorem 3.5.5]{CTS21}, 
the class $[X]\in \Br(V)$ is uniquely determined by
the class $[X_K]\in \Br(K)$, so it is enough to compute the $K$-conic $X_K$.


\ble
The structure group $\PGL(2)_K$
of the generic fibre $\T_K$ of $\T\to V$ can be reduced to 
$T=R^1_{L/K}(\G_{m,L})$, where $L=K(\sqrt{t^2-4})$
and $T$ is embedded into $\PGL(2)_K$ as the centraliser of an element of $\SL(2)_K$
of trace $t\in K$.
\ele
{\em Proof.} Let $g\in\SL(2)(K)$ be a matrix of trace $t\in K$.
The field extension $L/K$ is the splitting field of the characteristic polynomial of $g$, and thus
the centraliser of $g$ in $\SL(2)_K$ is isomorphic to $T$.
As above, $T\cong T/\mu_2$ is also the centraliser of $g$ in $\PGL(2)_K$, giving 
an embedding $T\hookrightarrow\PGL(2)_K$. 

The projection to the second factor $\SL(2)\times_k\SL(2)\to\SL(2)$ sends the generic
fibre $\T_K$ to the $\PGL(2)_K$-orbit of $g\in \SL(2)(K)$, where $\PGL(2)_K$
acts by conjugation. Let $\T_1$ be the preimage of $g$ in $\T_K$.
The $\PGL(2)_K$-orbit of $g$ is $\PGL(2)_K$-equivariantly isomorphic to
$T\backslash\PGL(2)_K$ so that $g$ is identified with the trivial coset of $T$.
 By Proposition \ref{comm} and its proof, 
$\T_1$ is a $K$-torsor for $T$ such that
$\T_K$ is the push-forward of $\T_1$ along $T\hookrightarrow\PGL(2)_K$. \hfill $\Box$


\medskip

By virtue of this lemma and Example \ref{dim2} it remains to compute $\T_1$. 
We shall do a little bit more
and describe an affine conic bundle $\pi\colon Y\to \A^3_k$ with generic fibre $Y_K=\T_1$.
The trace map $\tr\colon\SL(2)\to\A^1_k$ has a section $\si\colon \A^1_k\to \SL(2)$
that sends $t\in \bar k$ to the companion matrix 
$$g_t=\left(\begin{array}{cr}0&-1\\1&t\end{array}\right).$$ 
Let $Y$ be the preimage of $\si(\A^1_k)$
under the second projection $\SL(2)\times_k\SL(2)\to\SL(2)$.
Thus the $\bar k$-points of $Y$ are pairs of 
$(2\times 2)$-matrices $(A,g_t)$, where $A\in \SL(2)(\bar k)$ and $t\in \bar k$. 
As above, every fibre $\pi^{-1}(s,t,u)$ is a torsor for the centraliser of $g_t$ in $\PGL(2)$.
This centraliser is a torus if $g_t$ is semisimple, otherwise it is isomorphic to $\G_{a,k}$.
In particular, the generic fibre $Y_K=\T_1$ is a $K$-torsor for the torus $T=R^1_{L/K}(\G_{m,L})$.

\bpr \label{p1}
All fibres of $\pi\colon Y\to \A^3_k$ are non-empty affine conics, in particular, $\pi$ is flat and surjective.
The fibres of $\pi$ above the points of $V$
are affine subsets of smooth projective conics, the fibres above the smooth points of $M_0$
are singular conics consisting of two $\bar k$-lines meeting at a point, and the fibres
above the singular points of $M_0$ are double lines. The restriction of $\pi$ to the 
open subset of $\A^3_k$ given by $t^2-4\neq 0$ is isomorphic to 
\begin{equation}
\mu^2-(t^2-4)\nu^2=F(s,t,u).\label{equ1}
\end{equation}
The generic fibre $Y_K=\T_1$ of $\pi$ is a $K$-torsor for $T=R^1_{L/K}(\G_{m,L})$ 
isomorphic to the affine conic given by $(\ref{equ1})$.
\epr
{\em Proof.} The fibre $\pi^{-1}\big((s,t,u)\big)$ consists of the pairs $(M,g_t)$, where
$$M=\left(\begin{array}{cc}x+s&-y\\-tx-y-u&-x\end{array}\right)$$
is subject to $\det(M)=-(x^2+y^2+txy+sx+uy)=1$. This shows that $\pi$ is surjective
with all  fibres of dimension 1, hence $\pi$ is flat.
Homogenising this equation, we get the quadratic form
$$x^2+y^2+z^2+txy+sxz+uyz=0$$ whose discriminant is $-\frac{1}{4}F(s,t,u)$.
Thus the fibres of $\pi$ above the points of $V=\A^3_k\setminus M_0$ are affine subsets
of smooth projective conics.
Diagonalising this quadratic form, we set
$$ \mu=-\frac{t^2-4}{2}y+\frac{2u-st}{2}z, \quad \nu=-x-\frac{t}{2}y-\frac{s}{2}z,\quad \xi=z,$$
and obtain
$\mu^2-(t^2-4)\nu^2-F(s,t,u)\xi^2=0$.
This linear change of variables is given by a matrix with determinant $-\frac{1}{2}(t^2-4)$.
Setting $\xi=z=1$ we obtain (\ref{equ1}).
\hfill $\Box$

\bco \label{c1}
The class of the generic fibre of the smooth conic bundle $X\to V$ is
$$[X_K]=(t^2-4,F(s,t,u))\in\H^2(K,\mu_2).$$
\eco
{\em Proof.} This is immediate from Proposition \ref{p1}. Indeed,
as explained in Example \ref{dim2}, the conic $X_K$ is given by
the homogenisation of equation (\ref{equ1}).
\hfill $\Box$

\medskip

Note that the equation of the Markoff surface $M_d$ can be written as
\begin{equation}
(s^2-4)(t^2-4)=(2u-st)^2-4d. \label{norm1}
\end{equation}
For $d=0$ we obtain that the class of $F(s,t,u)$ in $k[s,t,u]/(t\pm 2)$ is a square, so $[X_K]$
is indeed unramified on $\A^3_k$ away from $M_0$ (which of course follows from 
the construction of $X$ from the $V$-torsor $\T$). It turns out that $[X_K]$
is also unramified at infinity.

\bpr
Let $M'_0\subset\P^3_k$ be the Zariski closure of $M_0$, and let $V'=\P^3_k\setminus M'_0$.
Then $[X]\in\Br(V)$ is contained in $\Br(V')\subset\Br(V)$.
\epr
{\em Proof.} Let $F'(r,s,t,u)$ be a homogeneous cubic polynomial such that $F(s,t,u)=F'(1,s,t,u)$. Rewriting the expression for $[X_K]$ from Corollary \ref{c1} in homogeneous
coordinates, that is, replacing $s,t,u$ by $s/r, t/r, u/r$, respectively, we obtain
$$[X_K]=\left(\frac{t^2-4r^2}{r^2},\frac{F'(r,s,t,u)}{r^3}\right)=
\left(\left(1-4\rho^2\right)\rho^{-2},\frac{F'(r,s,t,u)}{t^3}\rho^{-3}\right),$$
where $\rho=r/t$.
The rational function $\rho\in k(\P^3_k)$ is a uniformiser of the local ring of 
$\P^3_k$ at the divisor at infinity $r=0$, whereas $F'(r,s,t,u)t^{-3}\in k(\P^3_k)$ 
is a unit of this ring because $F'(0,s,t,u)=-stu\neq 0$.
We have $[X_K]=(1-4\rho^2,F'(r,s,t,u)t^{-3}\rho)$, hence
the residue of $[X_K]$ at $r=0$ is trivial. 
By the purity theorem for the Brauer group \cite[Theorem 3.7.1, formula (3.16)]{CTS21}, this implies the statement. \hfill $\Box$

\medskip

It would be interesting to find a conceptual explanation for the fact that the Brauer class $[X]$
is unramified at infinity.

\section{Word equations}
\subsection{Main theorem}

Let $w(x,y)$ be a non-trivial word in two letters, 
i.e.~a non-trivial element of the free group $\sF_2$ with two generators.
We write $w\colon \SL(2)\times\SL(2)\to\SL(2)$ for the $\PGL(2)$-equivariant
map defined by $w(x,y)$.
By a general theorem of Borel \cite[Theorem B]{Bor}, the morphism $w$ is dominant.

The Cayley--Hamilton theorem (\ref{CH}) gives rise to
a polynomial $P_w\in \Z[s,t,u]$ such that $\tr(w(A,B))=P_w(s,t,u)$, where
$s=\tr(A)$, $t=\tr(B)$, $u=\tr(AB)$. In other words,
the following diagram commutes:
\begin{equation}\begin{split}
\xymatrix{\SL(2)\times\SL(2)\ar[r]^{\ \ \ \ w}\ar[d]_{f}&\SL(2)\ar[d]^{\tr}\\
\A^3_k\ar[r]^{P_w}&\A^1_k} \label{d1}
\end{split}\end{equation}
The composition $\tr\circ w$ is dominant, so the polynomial $P_w(s,t,u)$ is non-constant.
For example, if $w(x,y)=xyx^{-1}y^{-1}$, then $P_w(s,t,u)=F(s,t,u)+2$. Moreover,
we have the following lemma.

\ble\label{divisible}
If $w\in [\sF_2,\sF_2]$, then there is a polynomial $Q_w(s,t,u)\in k[s,t,u]$ such that we have
\begin{equation} \label{decomp}
P_w(s,t,u)=Q_w(s,t,u)F(s,t,u)+2.
\end{equation}
\ele
{\em Proof.} 
Take any point $P=(s,t,u)\in M_0(\bar k)$. Let $\lambda,\mu\in \bar k$ be such that 
$\lambda^2-s\lambda+1=0$ and $\mu^2-t\mu+1=0$.
The Fricke identity (\ref{Fricke}) implies that $f^{-1}(P)$ contains
$$A=\left(\begin{matrix} \lambda& 0\\0&\lambda^{-1}\end{matrix}\right), \quad
B=\left(\begin{matrix} \mu& 0\\0&\mu^{-1}\end{matrix}\right).$$
Since $w(x,y)\in [\sF_2,\sF_2]$ and $AB=BA$, we have $w(A,B)=I$ so that $P_w(s,t,u)=\tr(I)=2$.
Thus, the restriction of $P_w(s,t,u)-2$ to the Cayley cubic $M_0$ is zero.
Since $F(s,t,u)$ is irreducible, we obtain \eqref{decomp}. \hfill $\Box$

\medskip

Let $\alpha\in \SL(2)(k)$ be a non-central element. 
Let $S_{w,\alpha}\subset\SL(2)\times\SL(2)$ be the closed subset given by
$w(A,B)=\alpha$. The isomorphism class of $S_{w,\alpha}$ 
depends only on the similarity class of $\alpha$, and so it is determined by $\tr(\alpha)$.

For $a\in k$ let $H_{w,a}\subset\A^3_k$ be the affine surface given by $P_w(s,t,u)=a$.
Thus we have a natural morphism $f\colon S_{w,\alpha}\to H_{w,a}$ with $a=\tr(\alpha)$.
As was pointed out in \cite[Proposition 2.2]{BZ}, the variety $S_{w,\alpha}$ is non-empty
when $\alpha$ is non-central and semisimple. Indeed, $P_w$ is non-constant, hence
surjective on $\bar k$-points. Thus the composition $P_w\circ f$ is surjective.
Hence $S_{w,\alpha'}$ is non-empty for some $\alpha'\in \SL(2)$ with $\tr(\alpha')=\tr(\alpha)$.
If $\alpha$ is semisimple and non-central, then $\alpha'$ and $\alpha$ are conjugate in 
$\SL(2)(\bar k)$, so that $S_{w,\alpha}$ is non-empty.

\bthe \label{t1}
Let $k$ be a field of characteristic different from $2$. 
Let $w(x,y)$ be a non-trivial word in two letters. Let
$\alpha\in \SL(2)(k)$ be a non-central semisimple element with trace $a=\tr(\alpha)$.
Let $X_{w,a}\to H_{w,a}\cap V$ be the restriction of the smooth conic bundle 
$X\to V$  to $H_{w,a}\cap V$. Then we have the following statements.

{\rm (i)} There is an open embedding $S_{w,\alpha}\cap f^{-1}(V)\to X_{w,a}$
of schemes over $H_{w,a}\cap V$ with dense image.

{\rm (ii)} Assume that $w(x,y)$ is in the commutator subgroup
$[\sF_2,\sF_2]\subset\sF_2$. Then $H_{w,a}\subset V$ so that there is 
an open embedding $S_{w,\alpha}\to X_{w,a}$ of $H_{w,a}$-schemes with dense image.
In particular, we have $\dim(S_{w,\alpha})=3$.
If, moreover, the polynomial $P_w(s,t,u)-a$ is geometrically irreducible, 
then $S_{w,a}$ is geometrically integral.
\ethe
{\em Proof.} (i) The centraliser of $\alpha$ in $\PGL(2)$ is
the norm 1 torus $T=R^1_{k'/k}(\G_{m,k'})$, where $k'=k(\sqrt{a^2-4})$.


Write $H'_{w,a}=H_{w,a}\cap V$ and 
let $\T_{w,a}=f^{-1}(H'_{w,a})$ be the restriction 
to $H'_{w,a}$ of the $V$-torsor $f\colon\T\to V$ for $\PGL(2)$ defined in Lemma \ref{2.1}.
We claim that the structure group of $\T_{w,a}\to H'_{w,a}$ reduces to $T$. 

Since $\alpha$ is semisimple, all elements of $\SL(2)(\bar k)$ of trace $a$
are conjugate to $\alpha$. Thus
the $\PGL(2)$-equivariant morphism $w\colon \SL(2)\times\SL(2)\to\SL(2)$ 
maps $\T_{w,a}$ to the $\PGL(2)$-orbit of
$\alpha$ in $\SL(2)$, where $\PGL(2)$ acts by conjugation. 
There is a unique $\PGL(2)$-equivariant isomorphism of this orbit and $T\backslash\PGL(2)$
that sends $\alpha$ to the trivial coset $x_0$ of $T$. We obtain
a $\PGL(2)$-equivariant morphism 
$$\T_{w,a}\to (T\backslash\PGL(2))\times_k H'_{w,a}$$ of
schemes over $H'_{w,a}$, so we can apply Proposition \ref{comm} (using that $T$ is commutative).
The desired $H'_{w,a}$-torsor for $T$ lifting $\T_{w,a}$ is the preimage of 
$x_0\times_k H'_{w,a}$ in $\T_{w,a}$, but
this is exactly $S_{w,\alpha}\cap f^{-1}(V)$.
Now (i) follows from Lemma \ref{lem4}.

(ii) Since $a\neq 2$, from Lemma \ref{divisible} we obtain that $M_0\cap H_{w,a}=\emptyset$,
 so that $H_{w,a}\subset V$ and $H'_{w,a}=H_{w,a}$, hence
$S_{w,\alpha}\cap f^{-1}(V)= S_{w,\alpha}$. Since $H_{w,a}\subset\A^3_k$ 
is the zero set of a non-constant polynomial, every irreducible component
of $H_{w,a}$ has dimension 2. Thus $\dim(S_{w,\alpha})=\dim(X_{w,a})=3$.

When $P_w(s,t,u)-a$ is geometrically irreducible, the surface
$H_{w,a}$ is geometrically integral. Since $X_{w,a}$ is a smooth conic bundle over $H_{w,a}$,
we see that $X_{w,a}$ is geometrically integral, so $S_{w,\alpha}$ is geometrically integral too.
\hfill $\Box$

\medskip

To illustrate the practical aspect of Theorem \ref{t1} (ii) we state the following

\bco
Let $k$ be a field of characteristic different from $2$. 
Let $w(x,y)$ be a non-trivial word in two letters contained in the commutator subgroup
$[\sF_2,\sF_2]\subset\sF_2$. Suppose that $(s,t,u)\in V(k)$ is such that $t^2-4$ is a
non-zero square in $k$. Then for any $\alpha\in\SL(2)(k)$ with $\tr(\alpha)=P_w(s,t,u)$
the equation $w(x,y)=\alpha$ has a solution $(A,B)\in\SL(2)(k)\times\SL(2)(k)$
such that $\tr(A)=s$, $\tr(B)=t$, $\tr(AB)=u$.
\eco
{\em Proof.} By Corollary \ref{c1}, the fibre of $X\to V$ over the $k$-point $(s,t,u)$ is isomorphic to 
$\P^1_k$. Now Theorem \ref{t1} (ii) implies that the fibre of $f\colon S_{w,\alpha}\to\A^3_k$
over $(s,t,u)$ is isomorphic to $\P^1_k$ with two $k$-points removed, so it has a $k$-point.
\hfill $\Box$

\subsection{Brauer groups of projective Markoff surfaces} \label{S 3.2}

For $d\neq 0$, $d\neq -4$
the Markoff surface $M_d$ is a dense open subset of the smooth cubic surface 
$M'_d\subset\P^3_k$ given by
$$ r(s^2+t^2+u^2)-stu-(d+4)r^3=0.$$
The Brauer group of $M'_d$ was computed in \cite[Proposition 3.2]{CTWX} 
and \cite[Lemmas 3.1, 3.2]{LM}. Let us recall this computation. 
Projection from the line $r=t=0$ contained in $M'_d$ to $\P^1_k$
with coordinates $(r:t)$ is a conic bundle $\phi\colon M'_d\to\P^1_k$
with five geometric singular fibres. 
The $k$-fibre  $\phi^{-1}(\infty)$ at the point at infinity $\infty=(0:1)$ is the union of two $k$-lines 
given by $r=s=0$ and $r=u=0$.
The remaining singular fibres are above the $k$-points
$t/r=\pm 2$ and above $(t/r)^2=d+4$, the latter being either a degree 2 closed point or a union of two
$k$-points. For these fibres, the quadratic extension over which the components 
are defined is given by adjoining $\sqrt{d}$. The element
$$\sA_0:=((t/r)^2-4,d)=(1-4(r/t)^2,d)\in \Br(k(\P^1_k))$$
has residue $d$ at the points $t/r=\pm 2$ and is unramified at all other points of $\P^1_k$.
We also note that the value of $\sA_0$ at the point at infinity $(r:t)=(0:1)$ is trivial.
Let $\sA$ be the image of $\sA_0$ in $\Br(k(M'_d))$.
We deduce that $\sA$ is unramified on $M'_d$ and hence $\sA\in\Br(M'_d)$, see
\cite[Proposition 6.2.7]{CTS21}.
Moreover, the value of $\sA$ at any $k$-point of $\phi^{-1}(\infty)$ is trivial.

For `general' $d$, namely, when
$[k(\sqrt{d},\sqrt{d+4}):k]=4$, we have $\Br(M'_d)/\Br(k)\cong\Z/2$
with generator $\sA$. If $d+4$ is a square in $k$, but $d$ is not, then 
$\Br(M'_d)/\Br(k)\cong(\Z/2)^2$ is generated by $\sA$ and another element.
Finally, if $d$ or $d(d+4)$ is a square in $k$, then $\Br(M'_d)/\Br(k)=0$.

\subsection{Commutator word variety}

In this section $w(x,y)$ is the commutator word $xyx^{-1}y^{-1}$.
For $\alpha\in\SL(2)(k)$ we denote the variety 
$S_{w,\alpha}\subset\SL(2)\times\SL(2)$ defined by $w(A,B)=\alpha$ simply by $S_\alpha$.

\bthe \label{3.3}
Let $k$ be a field of characteristic different from $2$. 
Let $\alpha\in\SL(2)(k)$ be a non-central semisimple element, and let $d=\tr(\alpha)-2$. 
Let $X_d\to M_d$ be the restriction of the smooth conic bundle 
$X\to V$ to $M_d\subset V$, $d\neq 0$, $d\neq -4$. Then we have the following statements.

{\rm (i)} The class $[X_d]\in\Br(M_d)$ equals $\sA\in\Br(M'_d)\subset\Br(M_d)$.

{\rm (ii)} $S_\alpha$ is isomorphic to a dense open subset of $X_d$.

{\rm (iii)} $S_\alpha$ is geometrically integral.

{\rm (iv)} $S_{\alpha}$ is $k$-unirational. If $k$ is infinite, 
$S_{\alpha}$ has a Zariski dense set of $k$-points.

{\rm (v)} $S_{\alpha}$ is $k$-rational if and only if the Markoff surface $M_d$ is $k$-rational.
This is the case if and only if $d$ or $d(d+4)$ is a square in $k$. 
\ethe
{\em Proof.} (i) The Markoff surface $M_d$ is given by $F(s,t,u)=d$, so 
for $d\neq 0$ we have $M_d\subset V$. Corollary \ref{c1} implies that $[X_d]=(t^2-4,d)=\sA$. 

(ii) For $w(x,y)=xyx^{-1}y^{-1}$ we have $P_w(s,t,u)=F(s,t,u)+2$.
The trace surface $H_{w,a}$, where $a=\tr(\alpha)=d+2$, is the zero set of 
$P_w(s,t,u)-a=F(s,t,u)-d$, hence $H_{w,a}=M_d$. 
Now (ii) follows from Theorem \ref{t1} (ii). 

(iii) This follows from Theorem \ref{t1} (ii) as $M_d$ is geometrically integral for any $d$.

(iv) To prove that $S_{\alpha}$ is $k$-unirational, 
by (ii) it is enough to prove that $X_d$ is $k$-unirational.
We start by recalling the well-known $k$-unirationality of $M'_d$.
For this it is enough to construct a rational curve $C\subset M'_d$ 
which is a double section of the conic bundle $\phi\colon M'_d\to\P^1_k$, because the pullback 
$D:=M'_d\times_{\P^1_k}C\to C$  is a conic bundle 
with a section, hence $D$ is a $k$-rational variety dominating $M'_d$.
We can take $C$ to be the line $r=t=0$. Indeed, the $\bar k$-fibres of $\phi$ are 
residual conics to this line, so $C$ is a double section of $\phi$.

The $k$-point $x_0=(0:0:0:1)$ is contained in $C\cap\phi^{-1}(\infty)$. As we have seen above,
this implies $\sA(x_0)=0$. Since $x_0$ is in $C(k)$, it lifts to a $k$-point $y_0$ on $D$. Thus
we have a dominant morphism $g\colon D\to M'_d$ of generic degree 2, 
where $D$ is a smooth, projective, and $k$-rational surface with a $k$-point $y_0$ such that
$g(y_0)=x_0$. The $k$-rationality of $D$ implies that the natural map $\Br(k)\to\Br(D)$
is an isomorphism. Thus $g^*\sA\in\Br(k)$. But $(g^*\sA)(y_0)=\sA(x_0)=0$, hence $g^*\sA=0$.
By (i) this implies that the pullback of the conic bundle $X_d\to M_d$ to 
$g^{-1}(M_d)\subset D$ has a section,
hence it is birationally equivalent to $D\times_k\P^1_k$, and therefore is $k$-rational.

(v) 
If $d$ or $d(d+4)$ is a square in $k$, then $M'_d$ is $k$-rational.  (In the first case
$M'_d$ contains two skew lines defined over $k$ and in the second case
$M'_d$ contains two skew lines conjugate over $k$ and
individually defined over $k(\sqrt{d+4})$, see \cite[Remark 3.3]{CTWX}.)
This implies $\Br(M'_d)/\Br(k)=0$. Since $\sA$ vanishes at a $k$-point of $M'_d$, we have
$\sA=0$. Thus $X_d$ is birationally equivalent to 
$M_d\times_k\P^1_k$, and hence is a $k$-rational variety.
(If $d(d+4)$ is a square in $k$, 
then the $k$-rationality of $X_d$ is proved in \cite[Lemma 4]{RBKC} by an explicit computation.)

Finally, assume that neither $d$ nor $d(d+4)$ is a square in $k$. 
Define $L=k(\sqrt{d+4})$ if $d+4$ is not a square in $k$, otherwise let $L=k$.
We note that $d+4$ is a square in $L$, but $d$ is not, and this implies that
$\Br(M'_{d,L})/\Br(L)\simeq(\Z/2)^2$ is generated by the images of $\sA$ and some other element
$\sA_1\in \Br(M'_{d,L})$.
By Lichtenbaum's theorem, the kernel of 
the restriction map $\Br(L(M_d))\to \Br(L(X_d))$ is generated by $\sA$, see 
\cite[Proposition 7.1.3]{CTS21} and the exact sequence (7.3) of {\em loc.~cit.}
Thus the image of $\sA_1$ in $\Br(L(X_d))$ is 
unramified over $L$ and is non-zero modulo $\Br(L)$.
The unramified Brauer group is a birational invariant of smooth,
proper, regular varieties \cite[Corollary 6.2.11]{CTS21}, hence
any smooth projective model of $X_{d,L}$ is not $L$-rational. 
Thus $X_{d,L}$ is not $L$-rational, and this implies that $X_d$ is not $k$-rational. 
Since $\Br(M'_d)/\Br(k)\neq 0$, the smooth projective surface $M'_d$ is not $k$-rational,
hence $M_d$ is not $k$-rational too.
\hfill $\Box$

\medskip

The existence of a $k$-point in $S_\alpha$ is a particular case of a general result due to 
R.C.~Thompson
\cite[Theorem 2]{Th}. When $k$ is a number field, our method can be used to prove 
a local-to-global statement for rational points on $S_\alpha$.
Variants of the following proof can be found in the existing literature, see \cite[Remark 6.4]{CTPS}.

Recall that a variety over a field $k$ is {\em split} if it contains an open geometrically integral 
$k$-subscheme, see \cite[Definition 10.1.3]{CTS21}.

\bpr \label{BMO}
Let $k$ be a number field. 
Let $\alpha\in\SL(2)(k)$ be a non-central semisimple element.
The Brauer--Manin obstruction is the only obstruction to
weak approximation on any smooth and proper variety birationally equivalent to
$S_\alpha$. Moreover, if we exclude the case when
$\tr(\alpha)+2$ is a square in $k$ but $\tr(\alpha)-2$ is not, 
then $S_\alpha$ satisfies weak approximation.
\epr
{\em Proof.} 
Using Hironaka's theorem, we can find a smooth, proper, geometrically integral variety $X'_d$ 
over $k$ that contains $X_d$ as an open subset. Moreover, we can choose $X'_d$ such
that there is a morphism $\varphi\colon X'_d\to M'_d$ extending $X_d\to M_d$. 

Let $x$ be a point of $M'_d$ of codimension 1. The irreducibility of $X'_d$ implies that
$\varphi^{-1}(x)$ is a curve. 
Thus the restriction of $\varphi$ to $\Spec(R)$, where $R$ is the local ring of $x$ in 
$M'_d$, is proper and flat (by miracle flatness, since all fibres have dimension 1). 
We claim that $\varphi^{-1}(x)$ is split.

By \cite[Proposition 10.1.12]{CTS21}, the splitness
of the fibre at $x$ does not depend on the choice of a regular integral scheme, proper and flat over
$\Spec(R)$, with the same generic fibre as $X_d\to M_d$. By Theorem \ref{3.3} (i)
this generic fibre is a smooth conic with class $\sA\in\Br(M'_d)\subset\Br(R)$.
Now \cite[Lemma 10.2.1]{CTS21} gives us such a scheme for which the fibre at $x$
is a smooth, hence split conic. (In the notation of {\em loc.~cit.} we have an equation of type I
because an equation of type II would imply that $\sA$ is ramified at $x\in\Spec(R)$
contrary to the fact that $\sA\in\Br(R)$.) This proves that $\varphi^{-1}(x)$ is split.

Let us prove the first statement. Write $\Omega$ for the set of places of $k$ and let $k_v$ be the completion of 
$k$ at $v\in\Omega$. Since $\sA\in\Br(M'_d)$ and $M'_d$ is proper,
there exists a finite subset $S_0\subset\Omega$ 
such that $\sA(P_v)=0\in\Br(k_v)$ for every $v\notin S_0$ and every $P_v\in M'_d(k_v)$,
see \cite[Proposition 13.3.1 (iii)]{CTS21}.

Let $S$ be a finite set of places of $k$ and let $P_v\in X'_d(k_v)$, $v\in S$, be 
local points coming from an adelic point $(P_v)_{v\in\Omega}$ in the Brauer--Manin set 
$X'_d({\bf A}_k)^\Br$.
We want to approximate $(P_v)_{v\in S}$ by a $k$-point. 
Without loss of generality we can assume that $S_0\subset S$.
After a small deformation we can arrange that $P_v\in X_d$, for $v\in S$.

By functoriality of the Brauer--Manin set
we have $(\varphi(P_v))_{v\in\Omega}\in M'_d({\bf A}_k)^\Br$.
Since $M'_d$ is birationally equivalent to a smooth projective surface
that is a conic bundle over $\P^1_k$
with four degenerate $\bar k$-fibres, by a theorem of Salberger and Colliot-Th\'el\`ene
\cite{CT90} the Brauer--Manin obstruction is the only obstruction to
weak approximation on $M'_d$. Thus there is a point $Q\in M'_d(k)$ arbitrarily close 
to each $\varphi(P_v)$ in the local topology of $k_v$, for $v\in S$. 
Since $\varphi(P_v)\in M_d(k_v)$ for $v\in S$, we can ensure that $Q\in M_d(k)$. Then
the fibre $\varphi^{-1}(Q)$ is a smooth projective conic with class $\sA(Q)\in\Br(k)$.
Since $\varphi(P_v)$ is close to $Q$, this conic has a $k_v$-point for all $v\in S$. 
Since $S_0\subset S$, 
it also has $k_v$-points for all $v\notin S$. By the Minkowski--Hasse theorem, 
$\varphi^{-1}(Q)$ has a
$k$-point. Then $\varphi^{-1}(Q)\simeq\P^1_k$, so we can find a $k$-point on $\varphi^{-1}(Q)$
arbitrarily close to $P_v$ for $v\in S$. 

Let us prove the second statement. 
By the results recalled in Section \ref{S 3.2}, the conditions guarantee that $\Br(M'_d)/\Br(k)$
is generated by the image of $\sA$. By Theorem \ref{3.3} (i)
the generic fibre of $\varphi$ is a conic with class $[\sA]$, hence $[\sA]$ is contained in
the kernel of the natural map
$\varphi^*\colon\Br(M'_d)\to \Br(X'_d)$, see \cite[Proposition 7.1.3]{CTS21}.
The image of $\Br(X'_d)$ in the Brauer group of the generic fibre of $\varphi$ comes from 
the Brauer group of the function field of $M'_d$, since this generic fibre is a conic. This
says that the Brauer group $\Br(X'_d)$ is vertical \cite[Definition 11.1.1]{CTS21}.
As we saw in the beginning of the proof, the fibres of $\varphi\colon X'_d\to M'_d$ 
over the codimension 1 points of $M'_d$ are split. This implies that 
$\varphi^*\colon\Br(M'_d)\to \Br(X'_d)$ is surjective \cite[Corollary 11.1.6]{CTS21}.
We conclude that $\Br(X'_d)=\Br(k)$, so 
there is no Brauer--Manin obstruction on $X'_d$ in this case. \hfill $\Box$

\brem{\rm
A well-known consequence of Proposition \ref{BMO} is that 
smooth and proper varieties birationally equivalent to $S_\alpha$ satisfy {\em weak weak
approximation}, that is, weak approximation outside of a finite set of primes.}
\erem

\subsection{Complements}

Part (i) of the following proposition can be compared to
a theorem of R.C. Thompson that $S_{-I}(k)\neq\emptyset$
if and only if $-1$ is a sum of two squares in $k$, see \cite[Theorem~1]{Th}.
Next, \cite[Theorem 2]{Th} says that
if $k$ has more than three elements, then every non-central element of
$\SL(2)(k)$ is a commutator. 
Part (ii) of the following proposition, together with Theorem \ref{3.3} (iv), gives a proof
of this result in the case of an infinite field $k$ of characteristic different from 2.

\bpr \label{compl}
Let $k$ be a field of characteristic different from $2$. 

{\rm (i)} The variety $S_{-I}$ is $k$-rational of dimension $3$
if $-1$ is a sum of two squares in $k$,
otherwise $S_{-I}(k)=\emptyset$ so that $S_{-I}$ is not $k$-rational.

{\rm (ii)} If $\alpha\in\SL(2)(k)$ is a non-central element
such that $\tr(\alpha)=\pm 2$, then $S_\alpha$ is $k$-rational of dimension $3$.
\epr
{\em Proof.} (i) If $A,B\in\SL(2)(\bar k)$ are such that
$ABA^{-1}B^{-1}=-I$, then $\tr(B)=\tr(-B)$, hence $\tr(B)=0$.
Similarly, we obtain $\tr(A)=0$. The Fricke identity (\ref{Fricke}) implies that $\tr(AB)=0$.
Thus $f(S_{-I})$ is a subset of $\{(0,0,0)\}$. By Lemma \ref{2.1}, 
$f^{-1}\big((0,0,0)\big)$ is a $k$-torsor for $\PGL(2)$. In particular,
it is not empty with transitive action of $\PGL(2)$, thus $f^{-1}\big((0,0,0)\big)=S_{-I}$.

By Corollary \ref{c1} the class of this torsor in the Brauer group $\Br(k)$
is $(-4,-4)=(-1,-1)$.
By a basic property of quaternion algebras, $(-1,-1)=0$ if and only if $-1$ is a norm of the quadratic extension $k(\sqrt{-1})$, that is, if and only if $-1$ a sum of two squares in $k$, 
cf.~\cite[Proposition 1.1.8]{CTS21}.

If the $k$-torsor $S_{-I}$ is trivial, it is isomorphic to $\PGL(2)$, in particular, it is $k$-rational. 
In the opposite case, it has no $k$-points, and hence is not $k$-rational.

(ii) The case $\tr(\alpha)=2$, $\alpha\neq I$, can be dealt with by an explicit computation. 
We can find a basis in which
$$\alpha=\left(\begin{array}{cc}1&1\\0&1\end{array}\right).$$
Suppose that $A, B\in\SL(2)(k)$ are $(2\times 2)$-matrices such that $ABA^{-1}B^{-1}=\alpha$.
We have $\tr(B)=\tr(ABA^{-1})=\tr(\alpha B)$, and this implies that $B$ is upper-triangular.
Likewise, we obtain that $A$ is also upper-triangular. Writing
$$A=\left(\begin{array}{cc}\lambda&x\\0&\lambda^{-1}\end{array}\right), \quad
B=\left(\begin{array}{cc}\mu&y\\0&\mu^{-1}\end{array}\right),$$
we obtain that the matrix equation $ABA^{-1}B^{-1}=\alpha$ is equivalent to
$$x\lambda(1-\mu^2)-y\mu(1-\lambda^2)=1.$$
We conclude that $S_\alpha$ is birationally equivalent to $\A^3_k$.

If $\tr(\alpha)=-2$, $\alpha\neq -I$, then our method still works. 
Let us explain how to adjust the proof of Theorem \ref{t1} to this case.
We have $a=\tr(\alpha)=d+2$ which gives $d=-4$, so $H_{w,a}=M_{-4}\subset V$. 
Let $M_\sharp:=M_{-4}\setminus\{(0,0,0)\}$. Let $\T_\sharp\to M_\sharp$
be the restriction of $\T\to V$ and let $X_\sharp\to M_\sharp$ be the restriction of $X\to V$.

The commutator gives a 
$\PGL(2)$-equivariant map $w\colon \SL(2)\times\SL(2)\to\SL(2)$
sending $\T_\sharp$ to
the subvariety of $\SL(2)$ whose $\bar k$-points are matrices of trace $-2$.
All such matrices, except $-I$, are conjugate to $\alpha$. The preimage $w^{-1}(-I)=S_{-I}$
is $f^{-1}\big((0,0,0)\big)$, by the proof of (i). Thus $w(\T_\sharp)$ is contained
in the $\PGL(2)$-orbit of $\alpha$ 
in $\SL(2)$, where $\PGL(2)$ acts by conjugation. Let $G$ be the centraliser of $\alpha$ in
$\PGL(2)$. Thus $G\cong\G_{a,k}$, in particular, $G$ is commutative.
There is a unique $\PGL(2)$-equivariant isomorphism of this orbit of $\alpha$
and $G\backslash\PGL(2)$ that sends $\alpha$ to the trivial coset $x_0$ of $G$. We obtain
a $\PGL(2)$-equivariant morphism 
$$\T_\sharp\to (G\backslash\PGL(2))\times_k M_\sharp$$ of
schemes over $M_\sharp$, so we can apply Proposition \ref{comm} 
(using that $G$ is commutative).
The desired $M_\sharp$-torsor for $G$ lifting $\T_\sharp$ is the preimage of 
$x_0\times_k M_\sharp$ in $\T_\sharp$, but this is exactly $S_\alpha$ because $f(S_\alpha)$
does not contain $(0,0,0)$. We know that $S_\alpha$
 is a dense open subset of $X_\sharp$ by Remark \ref{uni}. Moreover,
by the same remark, the generic fibre of $X_{-4}\to M_{-4}$ is isomorphic to the projective
line over the function field of $M_{-4}$, so  $S_\alpha$ is birationally equivalent to
$M_{-4}\times_k\P^1_k$. Finally, $M_{-4}$ is $k$-rational as a cubic surface with
a double $k$-point. Thus $S_\alpha$ is birationally equivalent to $\A^3_k$.
\hfill $\Box$

\medskip

When $n$ is prime, for any field $k$, Larsen and Lu showed that the commutator morphism
$w\colon\SL(n)\times_k\SL(n)\to\SL(n)$ is flat over the complement to the identity in $\SL(n)$,
 see \cite{LL21}.
Without using \cite{LL21} we have the following by-product of our method.

\bco \label{c2}
Let $k$ be a field of characteristic different from $2$. 

{\rm (i)} Let $w(x,y)$ be a non-trivial word in two letters
contained in the commutator subgroup $[\sF_2,\sF_2]\subset\sF_2$. Then the restriction of
the morphism $$w\colon\SL(2)\times_k\SL(2)\to\SL(2)$$ to the open subset
of non-central semisimple elements in $\SL(2)$  is faithfully flat.

{\rm (ii)} The commutator map $\SL(2)\times_k\SL(2)\to\SL(2)$ is faithfully flat over 
the complement to the identity in $\SL(2)$.
\eco
{\em Proof.} Since the source and the target are smooth, by miracle flatness
it is enough to check that all fibres 
have the same dimension. Thus (i) holds by Theorem \ref{t1} (ii), whereas 
(ii) holds by Theorem \ref{3.3} and Proposition \ref{compl}. \hfill $\Box$

\medskip

It is well-known that $S_I$ has dimension 4, so the result of Corollary \ref{c2} (ii) is best possible.
For the sake of completeness we note that $S_I$ is $k$-rational.

\bpr
The variety $S_I$ of commuting pairs of elements of $\SL(2)$ is $k$-rational of dimension $4$.
\epr
{\em Proof.} 
Motzkin and Taussky proved that $S_I$ is geometrically irreducible
\cite[page 399]{MT55}. The field of functions
$k(S_I)$ is the field of functions of the centraliser of the generic point $\Spec(K)$ of $\SL(2)$
in $\SL(2)_K$. The generic point is a regular semisimple element, so its
centraliser is a maximal torus in $\SL(2)$ defined over $K$. 
On the one hand, $\SL(2)$ is $k$-rational of dimension $3$, so $K$ is a purely transcendental
extension of $k$.
On the other hand, any maximal torus of $\SL(2)$ defined over $K$
is $K$-rational of dimension 1. This finishes the proof. \hfill $\Box$

Department of
Mathematics, Bar-Ilan University, 5290002 Ramat Gan, Israel

\texttt{bandman@macs.biu.ac.il}

\medskip

Department of
Mathematics, Bar-Ilan University, 5290002 Ramat Gan, Israel

\texttt{kunyav@macs.biu.ac.il}

\medskip

Department of Mathematics,
South Kensington Campus, Imperial College London SW7 2AZ,
United Kingdom -- and -- Institute for the Information Transmission
Problems, Russian Academy of Sciences, Moscow 127994, Russia

\texttt{a.skorobogatov@imperial.ac.uk}


{\small
\begin{thebibliography}{Mum74}
\itemsep=0pt

\bibitem[BZ16]{BZ} T. Bandman and Yu.G. Zarhin.
Surjectivity of certain word maps on $PSL(2,\mathbb C)$ and $SL(2,\mathbb C)$. {\em
Eur. J. Math.} {\bf 2} (2016) 614--643.

\bibitem[Bor83]{Bor} A. Borel. On free subgroups of semisimple groups. {\em Enseign.
Math.} {\bf 29} (1983) 151--164; reproduced in {\OE}uvres - Collected
Papers, vol.~IV, Springer-Verlag, Berlin--Heidelberg, 2001, pp.~41--54.

\bibitem[CT90]{CT90} J.-L. Colliot-Th\'el\`ene. Surfaces rationnelles fibr\'ees en coniques de degr\'e~4. {\em S\'em. th\'eorie des nombres Paris} 1988--89. Progress in Math., Birkh\"auser, 1990, p. 43--55.

\bibitem[CTPS]{CTPS} J.-L. Colliot-Th\'el\`ene, A. Pál, and A.N. Skorobogatov. Pathologies of the Brauer--Manin obstruction. {\em Math. Z.} {\bf 282} (2016) 799--817.

\bibitem[CTS21]{CTS21} J.-L. Colliot-Th\'el\`ene and A.N. Skorobogatov.
{\it The Brauer--Grothendieck group}, Springer, Cham, 2021.

\bibitem[CTWX20]{CTWX} J.-L. Colliot-Thélène, D. Wei, and F. Xu.
Brauer--Manin obstruction for Markoff surfaces. {\em Ann. Sc. Norm. Super. Pisa Cl. Sci.} (5) 
{\bf 21} (2020) 1257--1313.

\bibitem[Fri83]{Fri83} S. Friedland. Simultaneous similarity of matrices.
{\em Adv.~Math.} {\bf 50} (1983) 189--265.

\bibitem[GMS21]{GMS21} A. Ghosh, C. Meiri, and P. Sarnak. Commutators in ${\rm SL}_2$ and Markoff surfaces I. {\em New Zealand J. Math.} {\bf 52} (2021 [2021--2022]) 773--819.

\bibitem[GS22]{GS22}  A. Ghosh and P. Sarnak. Integral points on Markoff type cubic surfaces. 
{\em Invent. Math.} {\bf 229} (2022) 689--749. 

\bibitem[Gir71]{Giraud} J. Giraud. {\it Cohomologie non ab\'elienne.} Springer-Verlag, 1971.

\bibitem[GKP18]{GKP}
N.L. Gordeev, B.\`E. Kunyavski\u\i, and E.B. Plotkin.
Geometry of word equations in simple algebraic groups over special fields.
{\em Uspekhi Mat. Nauk} {\bf 73} (2018) 3--52;
English translation: {\em Russian Math. Surveys} {\bf 73} (2018) 753--796.

\bibitem[LL21]{LL21} M. Larsen and Z. Lu. Flatness of the commutator map over 
${\rm SL}_n$. {\em Int. Math. Res. Notices} (2021) 5605--5622.

\bibitem[LM21]{LM} D. Loughran and V. Mitankin. Integral Hasse principle and 
strong approximation for Markoff surfaces. {\em Int. Math. Res. Notices} (2021), no. 18, 14086--14122.

\bibitem[MT55]{MT55} T.S. Motzkin and O. Taussky. Pairs of matrices with property L. II. 
{\em Trans. Amer. Math. Soc.} {\bf 80} (1955) 387--401.

\bibitem[N00]{N} K. Nakamoto. Representation varieties and character varieties. 
{\em Publ. Res. Inst. Math. Sci.} {\bf 36} (2000) 159--189.

\bibitem[RBKC96]{RBKC} A.S. Rapinchuk, V.V. Benyash-Krivetz, and V.I. Chernousov.
Representation varieties of the fundamental groups of compact
orientable surfaces. {\em Israel J. Math.} {\bf 93} (1996) 29--71.

\bibitem[She84]{Shemesh}  D. Shemesh. Common eigenvectors of two matrices. 
{\em Linear Algebra Appl.} {\bf 62} (1984) 11--18.

\bibitem[Sko01]{S01} A. Skorobogatov. {\it Torsors and rational points.} Cambridge University Press,
2001.

\bibitem[Tho61]{Th} R.C. Thompson. Commutators in the special and general linear groups.
{\em Trans. Amer. Math. Soc.} {\bf 101} (1961) 16--33.

\end{thebibliography}
}
\end{document}